\documentclass[conference]{IEEEtran}
\ifCLASSINFOpdf
\else
\fi

\usepackage{amsmath, amssymb}

\newtheorem{theorem}{Theorem}
\newtheorem{proposition}[theorem]{Proposition}

\def \R {\mathbb{R}}

\def \E {\mathbb{E}}
\def \P {\mathbb{P}}

\def \NN {\mathcal{N}}

\def \e {\varepsilon}

\def \d {\delta}

\def \< {\langle}
\def \> {\rangle}
\def \^ {\widehat}

\def \dist {{\rm dist}}

\def \Span {{\rm span}}

\def \supp {{\rm supp}}
\def \Sparse {{\mathit{Sparse}}}

\begin{document}
%
\title{On the Role of Sparsity in Compressed Sensing and Random Matrix Theory}

\author{\IEEEauthorblockN{Roman Vershynin}
\IEEEauthorblockA{Department of Mathematics\\
University of Michigan\\
Ann Arbor, MI 48109\\
Email: romanv@umich.edu}}

\maketitle

\begin{abstract}
 We discuss applications of some concepts of Compressed Sensing 
 in the recent work on invertibility of random matrices due to Rudelson and the author.
 We sketch an argument leading to the optimal bound $\Omega(N^{-1/2})$ on the median 
 of the smallest singular value of an $N \times N$ matrix with random independent 
 entries. We highlight the parts of the argument where sparsity ideas played a key role.  
\end{abstract}

\IEEEpeerreviewmaketitle

\section{Introduction}

A concept that underlies the many recent developments in the area of Compressed Sensing 
is {\em sparsity}. Much earlier, sparsity has been used in a similar way (but often implicitly) 
in theoretical mathematics, and most notably in Geometric Functional Analysis. Recently the 
understanding of the role of sparsity led to formalizing some of the connections between the 
statements in those areas, leading to a new interplay between ``pure'' and ``applied'' 
mathematics.

While applications of ``pure'' mathematics to Compressed Sensing are expected and indeed
quite common, the reverse direction -- from Compressed Sensing to mathematics -- is still
rarely seen. 
This paper discusses one such application to the problem of {\em invertibility of random 
matrices}, which has been addressed in particular in the papers of Rudelson and the 
author \cite{RV square}, \cite{RV rectangular} and \cite{V products}.
Since we would like to focus here on techniques rather than results, we will often make 
oversimplifying assumptions and state weaker forms of the results. 
For the same reason, we discuss very little of history of these results and related work. 
The interested reader is encouraged to look at the original papers cited above for the 
statements of complete results, and for bibliography.

We will denote positive absolute constants by $C,c,C_1,\ldots$; their values may change from 
line to line.

\section{Sparsity as entropy control}

One normally thinks of sparsity as a way to represent objects (vectors or functions) 
in a certain basis in an economical way -- 
such that only a small number of basis elements can be used to accurately represent each object.
In this discussion, we shall identify the basis with the 
canonical basis of $\R^N$, and our objects will be vectors in $\R^N$.
We then say that a vector in $\R^N$ is {\em $s$-sparse} if it has few non-zero coordinates:
$$
|\supp(x)| \le s \ll N.
$$
We shall denote the set of all such vectors in $\R^N$ by $\Sparse(N,s)$. This set clearly consists of the union of all $s$-dimensional subspaces of $\R^N$.

An efficient way to use sparsity is through control of the metric entropy
of the space $\Sparse(N,s)$. 
Recall that, given a subset $S$ of a metric space and a number $\e >0$, 
the {\em covering number} $N(S,\e)$ is the smallest cardinality of an $\e$-net of $S$,
i.e. the smallest number of $\e$-balls centered at points in $S$ needed to cover $S$.
The logarithm of the covering number is often called {\em metric entropy}.

A simple argument based on comparison of volumes leads to an exponential 
bound of the metric entropy of many natural subsets of $\R^N$, and in particular of the Euclidean 
sphere $S^{N-1}$, see e.g. Lemma~9.5 in \cite{LT}.
This bound for the interesting range $\e \in (0,1)$
reads as
\begin{equation}									\label{covering sphere}
  N(S^{N-1}, \e) \le (3/\e)^N.
\end{equation}
This bound improves significantly for the set of sparse vectors. Since there are 
$\binom{N}{s}$ ways to choose the support of a sparse vector, we have
$$
N \big( \Sparse(N,s) \cap S^{N-1}, \e \big) \le \binom{N}{s} N(S^{N - 1},\e).
$$
Using \eqref{covering sphere} along with the bound 
$\binom{N}{s} \le (eN/s)^s$ valid for $1 \le s \le N/2$, which follows from
Stirling's formula, we conclude with
\begin{equation}									\label{sparse net}
  N \big( \Sparse(N,s) \cap S^{N-1}, \e \big) \le (CN/s)^s.
\end{equation}
Comparing this with \eqref{covering sphere}, we see that sparse vectors
enjoy significantly smaller entropy than the whole sphere --
the covering number is essentially exponential in the sparsity $s$ rather than the dimension $N$. 
This advantage is crucially used in many arguments, such as in the following one.

\medskip

In Compressed Sensing, a basic quality of matrices that guarantees their good 
performance as measurement operators is the {\em Restricted Isometry Condition}.
An $n \times N$ matrix $A$ with $n \le N$ is said to satisfy the Restricted Isometry Condition (RIC)
if $A$ acts as an approximate isometry when restricted to the set of sparse vectors. 
Formally, for every integer $s \le N$ we define the RIC constant $\d_s$ of the matrix $A$ 
as the minimal number that satisfies the two-sided inequality
\begin{equation}				\label{RIC}
  (1-\d) \|x\|_2^2 \le \|Ax\|_2^2 \le (1+\d) \|x\|_2^2
\end{equation}
for all $x \in \Sparse(N,s)$.
Candes and Tao \cite{CT} have shown (with constant improved in \cite{C}) 
that, given a matrix $A$ with $\d_{2s} \le \sqrt{2}-1$, one can exactly 
recover every $s$-sparse vector $x$ from its ``measurement vector'' $y=Ax$ by solving 
the convex optimization problem 
$$
\min_{x \in \R^N} \|x\|_1 \quad \text{subject to} \quad Ax=y.
$$

While it is difficult to explicitly construct matrices with good dimensions and RIC 
constants, random constructions are abundant in the literature 
(see Section V in \cite{CW}). Here we sketch the 
known argument for Gaussian matrices, which will highlight 
sparsity as entropy control, and will lead to our discussion of more difficult 
questions in Random Matrix Theory.

\medskip

\begin{proposition}[Gaussian matrices]					\label{gaussian RIC}
  Let $\bar{A}$ be an $n \times N$ matrix whose entries are independent standard normal 
  random variables. Let $1 \le s \le N$ and $\d > 0$. If 
  $$
  n \ge C(\d) s \log (N/s)
  $$
  then, with high probability, the matrix 
  $A = \frac{1}{\sqrt{n}} \bar{A}$ satisfies RIC with constant $\d_s \le \d$.
  Here $C(\d)>0$ only depends on $\d$.
\end{proposition}

\smallskip

\begin{proof}(Sketch)
An approximation argument shows that it is enough to check \eqref{RIC} for 
all $x$ in any fixed $\d$-net of $\Sparse(N,s) \cap S^{N-1}$.
So we choose such a net $\NN$ of cardinality controlled as in \eqref{sparse net},
and we fix a vector $x \in \NN$.
Due to independence of the rows of $A$ and the rotation invariance of the normal 
distribution, the random variable $\|Ax\|_2^2$ is distributed identically
with $\chi^2 := \sum_{i=1}^n g_i^2$,
where $g_i$ are independent standard normal random variables. 
By the known concentration properties of the $\chi^2$ distribution, or alternatively 
by the standard exponential concentration inequalities, one has 
$$
(1-\d) n \le \sum_{i=1}^n g_i^2 \le (1+\d) n
$$
with probability at least $1 - e^{-c(\d) n}$.
In other words, with this probability, the Restricted Isometry Condition \eqref{RIC}
holds for a fixed vector $x \in \NN$. Taking the union bound and using 
the bound \eqref{sparse net} on the cardinality of the net, we see that \eqref{RIC}
holds for {\em all} vectors $x \in \NN$ with probability at least
$$
1 - |\NN| e^{-c(\d) n} \ge 1 - (CN/s)^s e^{-c(\d) n}.
$$
By the condition we made on the dimensions, the proof is complete.
\end{proof}

The argument above can be easily generalized to distributions other than normal
by using standard exponential concentration inequalities; 
suitable moment bounds (subgaussian) are sufficient for this purpose.

\section{Invertibility of random matrices}

One can view the Restricted Isometry Condition \eqref{RIC} 
as the condition that all submatrices of $A$ with a given 
number of columns are well conditioned. 
The question of how well conditioned random matrices are goes back to 
at least Von Neumann and his collaborators, in connection with their
work on large matrix inversion. Some history of the work on this 
problem is described in \cite{RV square} and \cite{RV rectangular}, 
and some new results appeared since then, see \cite{TV square universality}.
Here we shall focus on the original prediction going back to Von Neumann
and his group -- that the smallest singular value $s_N(A)$ of an 
$N \times N$ matrix with random independent centered entries is 
typically of order $N^{-1/2}$. Coupled with the known 
estimate on the largest singular value $\E s_1(A) \le N^{1/2}$ (valid
under suitable moment assumptions), the prediction implies
that the condition number $\kappa(A) = s_1(A) s_N(A) = O(N)$, 
i.e. is typically linear in the dimension.

This prediction was verified for Gaussian matrices in \cite{E, S}
using the explicit formula for the joint density of their eigenvalues, 
and was first proved for general random matrices in \cite{RV square}
under some mild moment assumptions. Ideas based on sparsity
play an important role in \cite{RV square}. We will discuss this role the in the 
rest of the paper, and sketch the proof of the prediction above:

\medskip

\begin{theorem}[\cite{RV square}]						\label{square}
  Let $A$ be an $N \times N$ matrix whose entries are
  independent identically distributed random variables
  with mean zero, unit variance, and fourth moment bounded by 
  a constant. 
  Then the median of $s_N(A)$ is bounded below by $cN^{-1/2}$.
\end{theorem}

\medskip

Note that the result is sharp -- it was proved in \cite{RV square above}
that the median of $s_N(A)$ is bounded above by $CN^{-1/2}$.

\section{Invertibility on sparse vectors}

Our plan is to first prove Theorem~\ref{square} for Gaussian matrices $A$
(whose all entries are standard normal random variables), and then 
to indicate how to modify the proof for general distributions.

The smallest singular value has the following convenient expression:
$$
s_N(A) = \min_{x \in S^{N-1}} \|Ax\|_2.
$$
Our goal is then to bound $\|Ax\|_2$ below uniformly for all unit vectors $x$.

We already know how to achieve this goal for all {\em sparse} vectors $x$. 
Indeed, by Proposition~\ref{gaussian RIC} the Gaussian matrices satisfy 
the Restricted Isometry Property. 
If we choose $n=N$ and $s = cN$ with sufficiently small absolute constant $c > 0$, 
Proposition~\ref{gaussian RIC} shows that, with high probability,
$$
\min_{x \in \Sparse(N,cN) \cap S^{N-1}} \|Ax\|_2 \ge c N^{1/2}.
$$
Note that this bound is much better than we need in Theorem~\ref{square} -- 
we would be happy with $cN^{-1/2}$ in the right hand side.

Now we need to handle the non-sparse vectors.

\section{Invertibility on spread vectors}

Our success with sparse vectors is due the fact that there are ``not too many''
of them. As we have seen by comparing \eqref{covering sphere}
to \eqref{sparse net}, the metric entropy of the set of sparse vectors is much smaller
than that of all vectors. Such a nice entropy control allowed us to handle all sparse
vectors by taking a union bound (in the proof of Proposition~\ref{gaussian RIC}) 
without paying too much price in the probability estimates.

Repeating a similar argument for non-sparse vectors is hopeless, as they lack a
nice entropy control. Instead, we could first try to identify the class of vectors which is
entirely opposite to the sparse vectors, and try to handle this class. 
These are {\em spread vectors} -- those vectors in $S^{N-1}$ 
whose all coordinates have the same order $N^{-1/2}$.
An advantage of spread vectors over sparse ones is that we know the magnitude
of all their coefficients. So we develop the following {\em geometric} argument
to prove the invertibility on the set of spread vectors. 

Let us begin with a qualitative argument. Suppose the matrix $A$ performs extremely
poor, and we have $s_N(A) = 0$; in other words, $A$ is a singular matrix. 
Therefore one of its columns $A_k$ of $A$ lies in the span $H_k = \Span(A_i)_{i \ne k}$ 
of the others.

This simple observation can be made into a quantitative argument, which will work very 
well with the spread vectors. Suppose $x = (x_1,\ldots,x_N) \in \R^N$ is a spread vector. 
Then, for every $k=1,\ldots,N$, we have 
\begin{align}					\label{Ax spread}
\|Ax\|_2 
  &\ge \dist(Ax, H_k) 
  = \dist \Big( \sum_{i=1}^N x_i A_i, H_k \Big)  \\ 
  &= \dist (x_k A_k, H_k)
  = |x_k| \cdot \dist (X_k, H_k) \nonumber\\
  &\ge c N^{-1/2} \; \dist(X_k, H_k). \nonumber
\end{align}
Since the left hand side does not depend on $x$, we have proved in particular 
(for $k=1$) that
$$
\min_{\text{Spread } x} \|Ax\|_2 \ge c N^{-1/2} \; \dist(X_1, H_1).
$$

It remains to estimate the distance between the random vector $X_1$ 
and the independent random hyperplane $H_1$.
Since $X_1$ is a Gaussian vector, it is easy to check (using the rotation invariance
of the Gaussian disstibution) that $\dist(X_1, H_1)$ is distributed
identically with the absolute value of a standard normal random variable $g$.
But the density of $g$ is bounded by the absolute constant $(2\pi)^{-1/2}$, which 
makes 
$$
\dist(X_1, H_1) = |g| = \Omega(1)
$$
with arbitrarily high constant probability (say, $0.999$).

We have thus shown that, with arbitrarily high constant probability,
$$
\min_{\text{Spread } x} \|Ax\|_2 \ge c N^{-1/2}.
$$
This is a desired uniform bound for the spread vectors.

\section{Bridging sparse and spread vectors}

There are of course many vectors that are neither sparse nor spread, 
but it will now be relatively easy to bridge these two classes. 

Consider all vectors in $S^{N-1}$ that are within a small absolute constant distance 
$c'>0$ from the set of sparse vectors $\Sparse(N,cN) \cap S^{N-1}$.
We shall call such vectors {\em compressible}, and the rest of the vectors on the
sphere are {\em incompressible}. The intuition, which again is coming from sparse recovery,
suggests that compressible vectors should behave similarly to sparse vectors, 
while incompressible vectors should be similar to spread vectors. 

Indeed, a trivial approximation argument extends our invertibility bound from sparse to 
compressible vectors (one just need to approximate a compressible vector by a sparse one 
and use that the error of this approximation $c'$ can only blow up by a factor 
$\|A\| = O(N^{1/2})$). So we have the desired bound
$$
\min_{\text{Compressible } x} \|Ax\|_2 \ge c N^{1/2}.
$$

For incompressible vectors, instead of an approximation argument (which won't work)
one makes the following simple observation:  every incompressible vector 
has $\Omega(N)$ coordinates of magnitude $\Omega(N^{-1/2})$.
This is a way how incompressible vectors are similar to spread ones. 

To complete the proof for incompressible vectors, we again use the geometric argument,
but stop just before the last estimate in \eqref{Ax spread}:
$$
\|Ax\|_2 \ge \max_k |x_k| \cdot \dist(X_k, H_k).
$$
As we already know, for each $k=1,\ldots,N$, the distance satisfies
$\dist(X_k, H_k) = \Omega(1)$ with arbitrarily large probability. 
Therefore, still with high probability, most of these distances 
(arbitrarily high constant proportion of them) are $\Omega(1)$. 
On the other hand, we also know that some fixed proportion of the 
coordinates $x_k$ are of magnitude $\Omega(N^{-1/2})$. Therefore, 
intersecting these two events, we see that for any incompressible vector $x$
there exists a coordinate $k$ that satisfies both bounds. This implies
that, with high probability,
$$
\min_{\text{Incompressible } x} \|Ax\|_2 \ge c N^{-1/2}.
$$
This is a desired bound which, along with the already proved estimate for 
compressible vectors, implies the final result: 
$$
s_N(A) = \min_{x \in S^{N-1}} \|Ax\|_2
\ge cN^{-1/2}.
$$

\section{Extensions and further remarks}

The above argument generalizes from Gaussian to general distributions. 
There are two places where we used rotation invariance of the Gaussian distribution.
 One such place was the use of Proposition~\ref{gaussian RIC} in the treatment 
of sparse vectors. As we already mentioned, Proposition~\ref{gaussian RIC} 
can be easily extended to more general distributions using 
the standard large deviation inequalities.

The other place where Gaussian distribution was used was in the argument 
for spread vectors. We argued there that the distance $\dist(X_1, H_1)$ between 
a random vector and a random independent hyperplane is $\Omega(1)$ 
with arbitrarily high probability.  For Gaussian distribution of the entries, 
this followed by a direct and easy computation.
For more general distributions, this estimate is still true, but it requires
more work. 

Let us condition on a realization of the hyperplane $H_1$, and let 
$a \in \R^N$ be a normal vector of $H_1$. Then clearly
$$
\dist(X_1, H_1) = \< a, X_1 \> .
$$
Writing this is coordinates for $a = (a_1,\ldots,a_N)$ and 
$X_1 = (\xi_1, \ldots, \xi_N)$ we see that 
$$
\dist(X_1, H_1) =  \Big| \sum_{i=1}^N a_i \xi_i \Big| =: |S|
$$
where $S$ is clearly a sum of independent random variables. 

Our goal is to show that $|S| = \Omega(1)$ with high probability.
One way to do this is to use a Central Limit Theorem (in the form of Berry-Esseen)
to approximate $S$ by a standard normal random variable $g$, for which we
already have the desired result. For the The Central Limit Theorem to work, 
one obviously needs that many coordinates of $a$ are not too small (for example, 
it will clearly not work if $a$ is $1$-sparse, as the sum $S$ will consist of just one term).
However, sparsity ideas can be again of help here. Running an argument similar
to the one above for compressible vectors, one can show that, with high probability,
the normal $a$ to the random hyperplane $H_1$ is incompressible. 
We then condition on such $H_1$, and the Central Limit Theorem works well:
\begin{equation}					\label{CLT}
  \big| \P( |S| < \e) - \P(|g| < \e) \big| = O( N^{-1/2} ).
\end{equation}
This proves the desired bound $\dist(X_1, H_1) = \Omega(1)$ with high probability
$1 - O(N^{-1/2})$.

\medskip

In \cite{RV rectangular, V products}, Theorem~\ref{square} 
was extended to {\em rectangular matrices} $N \times n$, where $N \ge n$. 
Under the same assumptions, the median of the smallest singular value $s_n(A)$ 
of such random matrices is bounded below by $c ( \sqrt{N} - \sqrt{n-1} )$, 
which is asymptotically optimal.
Note that for square matrices, where $N=n$, this bound equals $c N^{-1/2}$, 
which agrees with Theorem~\ref{square}.

\medskip

Under stronger moment assumptions on the entries (subgaussian), not 
only the median of the smallest singular value can be estimated, but also 
strong probability inequalities can be proved. 
For example, square matrices satisfy 
$$
\P( s_N(A) < \e N^{-1/2} ) \le C \e + e^{-cN}
$$
and rectangular matrices satisfy 
$$
\P \big( s_n(A) \le \e (\sqrt{N} - \sqrt{n-1}) \big) 
\le (C\e)^{N-n+1} + e^{-cN}
$$
for all $\e \ge 0$.
Proving such exponential inequalities is more difficult, because 
one can not afford a polynomial error in probability $O(N^{-1/2})$
which one necessarily
obtains when applying Central Limit Theorem in \eqref{CLT}.
Instead of using Central Limit Theorems, one develops a Littlewood-Offord Theory, 
whose probability estimates are fine-tuned to the additive structure
of the coefficients of $a$. Since the sparsity does not play a key role in 
these arguments, we will not discuss this direction here. 
The interested reader is encouraged to consult the papers \cite{RV rectangular, V products}.

\section*{Acknowledgment}
This research was partially supported by NSF FRG grant DMS 0918623.

\end{document}